%
%
%
\def\unredoffs{} \def\redoffs{\voffset=-.31truein\hoffset=-.59truein}
\def\speclscape{\special{ps: landscape}}
%
%
%
%
\newbox\leftpage \newdimen\fullhsize \newdimen\hstitle \newdimen\hsbody
\tolerance=1000\hfuzz=2pt
\catcode`\@=11 
%
\ifx\answ\bigans\message{(This will come out unreduced.}
\magnification=1200\unredoffs\baselineskip=16pt plus 2pt minus 1pt
\hsbody=\hsize \hstitle=\hsize 
\else\message{(This will be reduced.} \let\l@r=L
\magnification=1000\baselineskip=16pt plus 2pt minus 1pt \vsize=7truein
\redoffs \hstitle=8truein\hsbody=4.75truein\fullhsize=10truein\hsize=\hsbody
\output={\ifnum\pageno=0 
  \shipout\vbox{\speclscape{\hsize\fullhsize\makeheadline}
    \hbox to \fullhsize{\hfill\pagebody\hfill}}\advancepageno
  \else
  \almostshipout{\leftline{\vbox{\pagebody\makefootline}}}\advancepageno
  \fi}
\def\almostshipout#1{\if L\l@r \count1=1 \message{[\the\count0.\the\count1]}
      \global\setbox\leftpage=#1 \global\let\l@r=R
 \else \count1=2
  \shipout\vbox{\speclscape{\hsize\fullhsize\makeheadline}
      \hbox to\fullhsize{\box\leftpage\hfil#1}}  \global\let\l@r=L\fi}
\fi
%
\newcount\yearltd\yearltd=\year\advance\yearltd by -1900

\def\Title#1#2{\nopagenumbers\abstractfont\hsize=\hstitle\rightline{#1}%
\vskip 1in\centerline{\titlefont #2}\abstractfont\vskip .5in\pageno=0}
%

%
\def\Date#1{\vfill\leftline{#1}\tenpoint\supereject\global\hsize=\hsbody%
\footline={\hss\tenrm\folio\hss}}
%

\def\draftmode{\message{ DRAFTMODE }\def\draftdate{{\rm preliminary draft:
\number\month/\number\day/\number\yearltd\ \ \hourmin}}%
\headline={\hfil\draftdate}\writelabels\baselineskip=20pt plus 2pt minus 2pt
 {\count255=\time\divide\count255 by 60 \xdef\hourmin{\number\count255}
  \multiply\count255 by-60\advance\count255 by\time
  \xdef\hourmin{\hourmin:\ifnum\count255<10 0\fi\the\count255}}}
\def\nolabels{\def\wrlabeL##1{}\def\eqlabeL##1{}\def\reflabeL##1{}}
\def\writelabels{\def\wrlabeL##1{\leavevmode\vadjust{\rlap{\smash%
{\line{{\escapechar=` \hfill\rlap{\sevenrm\hskip.03in\string##1}}}}}}}%
\def\eqlabeL##1{{\escapechar-1\rlap{\sevenrm\hskip.05in\string##1}}}%
\def\reflabeL##1{\noexpand\llap{\noexpand\sevenrm\string\string\string##1}}}
\nolabels
%
\global\newcount\secno \global\secno=0
\global\newcount\meqno \global\meqno=1
\def\newsec#1{\global\advance\secno by1\message{(\the\secno. #1)}
\global\subsecno=0\eqnres@t\noindent{\bf\the\secno. #1}
\writetoca{{\secsym} {#1}}\par\nobreak\medskip\nobreak}
\def\eqnres@t{\xdef\secsym{\the\secno.}\global\meqno=1\bigbreak\bigskip}
\def\sequentialequations{\def\eqnres@t{\bigbreak}}\xdef\secsym{}
\global\newcount\subsecno \global\subsecno=0
\def\subsec#1{\global\advance\subsecno by1\message{(\secsym\the\subsecno. #1)}
\ifnum\lastpenalty>9000\else\bigbreak\fi
\noindent{\it\secsym\the\subsecno. #1}\writetoca{\string\quad
{\secsym\the\subsecno.} {#1}}\par\nobreak\medskip\nobreak}
\def\appendix#1#2{\global\meqno=1\global\subsecno=0\xdef\secsym{\hbox{#1.}}
\bigbreak\bigskip\noindent{\bf Appendix #1. #2}\message{(#1. #2)}
\writetoca{Appendix {#1.} {#2}}\par\nobreak\medskip\nobreak}
%
%
\def\eqnn#1{\xdef #1{(\secsym\the\meqno)}\writedef{#1\leftbracket#1}%
\global\advance\meqno by1\wrlabeL#1}
\def\eqna#1{\xdef #1##1{\hbox{$(\secsym\the\meqno##1)$}}
\writedef{#1\numbersign1\leftbracket#1{\numbersign1}}%
\global\advance\meqno by1\wrlabeL{#1$\{\}$}}
\def\eqn#1#2{\xdef #1{(\secsym\the\meqno)}\writedef{#1\leftbracket#1}%
\global\advance\meqno by1$$#2\eqno#1\eqlabeL#1$$}
%
\newskip\footskip\footskip14pt plus 1pt minus 1pt 
\def\footnotefont{\ninepoint}\def\f@t#1{\footnotefont #1\@foot}
\def\f@@t{\baselineskip\footskip\bgroup\footnotefont\aftergroup\@foot\let\next}
\setbox\strutbox=\hbox{\vrule height9.5pt depth4.5pt width0pt}
\global\newcount\ftno \global\ftno=0
\def\foot{\global\advance\ftno by1\footnote{$^{\the\ftno}$}}
%
\newwrite\ftfile
\def\footend{\def\foot{\global\advance\ftno by1\chardef\wfile=\ftfile
$^{\the\ftno}$\ifnum\ftno=1\immediate\openout\ftfile=foots.tmp\fi%
\immediate\write\ftfile{\noexpand\smallskip%
\noexpand\item{f\the\ftno:\ }\pctsign}\findarg}%
\def\footatend{\vfill\eject\immediate\closeout\ftfile{\parindent=20pt
\centerline{\bf Footnotes}\nobreak\bigskip\input foots.tmp }}}
\def\footatend{}
%
%
\global\newcount\refno \global\refno=1
\newwrite\rfile
%
\def\ref{\nref}
\def\nref#1{\xdef#1{[\the\refno]}\writedef{#1\leftbracket#1}%
\ifnum\refno=1\immediate\openout\rfile=refs.tmp\fi
\global\advance\refno by1\chardef\wfile=\rfile\immediate
\write\rfile{\noexpand\item{#1\ }\reflabeL{#1\hskip.31in}\pctsign}\findarg}
\def\findarg#1#{\begingroup\obeylines\newlinechar=`\^^M\pass@rg}
{\obeylines\gdef\pass@rg#1{\writ@line\relax #1^^M\hbox{}^^M}%
\gdef\writ@line#1^^M{\expandafter\toks0\expandafter{\striprel@x #1}%
\edef\next{\the\toks0}\ifx\next\em@rk\let\next=\endgroup\else\ifx\next\empty%
\else\immediate\write\wfile{\the\toks0}\fi\let\next=\writ@line\fi\next\relax}}
\def\striprel@x#1{} \def\em@rk{\hbox{}}
\def\lref{\begingroup\obeylines\lr@f}
\def\lr@f#1#2{\gdef#1{\ref#1{#2}}\endgroup\unskip}

\def\addref#1{\immediate\write\rfile{\noexpand\item{}#1}} 
\def\startrefs#1{\immediate\openout\rfile=refs.tmp\refno=#1}
\def\refs#1{\count255=1[\r@fs #1{\hbox{}}]}
\def\r@fs#1{\ifx\und@fined#1\message{reflabel \string#1 is undefined.}%
\nref#1{need to supply reference \string#1.}\fi%
\vphantom{\hphantom{#1}}\edef\next{#1}\ifx\next\em@rk\def\next{}%
\else\ifx\next#1\ifodd\count255\relax\xref#1\count255=0\fi%
\else#1\count255=1\fi\let\next=\r@fs\fi\next}
%

%
\newwrite\ffile\global\newcount\figno \global\figno=1
\def\fig{fig.~\the\figno\nfig}
\def\nfig#1{\xdef#1{fig.~\the\figno}%
\writedef{#1\leftbracket fig.\noexpand~\the\figno}%
\ifnum\figno=1\immediate\openout\ffile=figs.tmp\fi\chardef\wfile=\ffile%
\immediate\write\ffile{\noexpand\medskip\noexpand\item{Fig.\ \the\figno. }
\reflabeL{#1\hskip.55in}\pctsign}\global\advance\figno by1\findarg}
\def\xfig{\expandafter\xf@g}\def\xf@g fig.\penalty\@M\ {}
\def\figs#1{figs.~\f@gs #1{\hbox{}}}
\def\f@gs#1{\edef\next{#1}\ifx\next\em@rk\def\next{}\else
\ifx\next#1\xfig #1\else#1\fi\let\next=\f@gs\fi\next}
\newwrite\lfile
{\escapechar-1\xdef\pctsign{\string\%}\xdef\leftbracket{\string\{}
\xdef\rightbracket{\string\}}\xdef\numbersign{\string\#}}

\def\writestop{\def\writestoppt{\immediate\write\lfile{\string\pageno%
\the\pageno\string\startrefs\leftbracket\the\refno\rightbracket%
\string\def\string\secsym\leftbracket\secsym\rightbracket%
\string\secno\the\secno\string\meqno\the\meqno}\immediate\closeout\lfile}}
\def\writestoppt{}\def\writedef#1{}
\def\seclab#1{\xdef #1{\the\secno}\writedef{#1\leftbracket#1}\wrlabeL{#1=#1}}
\def\subseclab#1{\xdef #1{\secsym\the\subsecno}%
\writedef{#1\leftbracket#1}\wrlabeL{#1=#1}}
\newwrite\tfile \def\writetoca#1{}
\def\leaderfill{\leaders\hbox to 1em{\hss.\hss}\hfill}
\def\writetoc{\immediate\openout\tfile=toc.tmp
   \def\writetoca##1{{\edef\next{\write\tfile{\noindent ##1
   \string\leaderfill {\noexpand\number\pageno} \par}}\next}}}
%
%
%

\catcode`\@=12 
%
\edef\tfontsize{\ifx\answ\bigans scaled\magstep3\else scaled\magstep4\fi}
\font\titlerm=cmr10 \tfontsize \font\titlerms=cmr7 \tfontsize
\font\titlermss=cmr5 \tfontsize \font\titlei=cmmi10 \tfontsize
\font\titleis=cmmi7 \tfontsize \font\titleiss=cmmi5 \tfontsize
\font\titlesy=cmsy10 \tfontsize \font\titlesys=cmsy7 \tfontsize
\font\titlesyss=cmsy5 \tfontsize \font\titleit=cmti10 \tfontsize
\skewchar\titlei='177 \skewchar\titleis='177 \skewchar\titleiss='177
\skewchar\titlesy='60 \skewchar\titlesys='60 \skewchar\titlesyss='60
\def\titlefont{\def\rm{\fam0\titlerm}
\textfont0=\titlerm \scriptfont0=\titlerms \scriptscriptfont0=\titlermss
\textfont1=\titlei \scriptfont1=\titleis \scriptscriptfont1=\titleiss
\textfont2=\titlesy \scriptfont2=\titlesys \scriptscriptfont2=\titlesyss
\textfont\itfam=\titleit \def\it{\fam\itfam\titleit}\rm}
 \ifx\answ\bigans\else scaled\magstep1\fi
\ifx\answ\bigans\def\abstractfont{\tenpoint}\else
\font\abssl=cmsl10 scaled \magstep1
\font\absrm=cmr10 scaled\magstep1 \font\absrms=cmr7 scaled\magstep1
\font\absrmss=cmr5 scaled\magstep1 \font\absi=cmmi10 scaled\magstep1
\font\absis=cmmi7 scaled\magstep1 \font\absiss=cmmi5 scaled\magstep1
\font\abssy=cmsy10 scaled\magstep1 \font\abssys=cmsy7 scaled\magstep1
\font\abssyss=cmsy5 scaled\magstep1 \font\absbf=cmbx10 scaled\magstep1
\skewchar\absi='177 \skewchar\absis='177 \skewchar\absiss='177
\skewchar\abssy='60 \skewchar\abssys='60 \skewchar\abssyss='60
\def\abstractfont{\def\rm{\fam0\absrm}
\textfont0=\absrm \scriptfont0=\absrms \scriptscriptfont0=\absrmss
\textfont1=\absi \scriptfont1=\absis \scriptscriptfont1=\absiss
\textfont2=\abssy \scriptfont2=\abssys \scriptscriptfont2=\abssyss
\textfont\itfam=\bigit \def\it{\fam\itfam\bigit}\def\footnotefont{\tenpoint}%
\textfont\slfam=\abssl \def\sl{\fam\slfam\abssl}%
\textfont\bffam=\absbf \def\bf{\fam\bffam\absbf}\rm}\fi
\def\tenpoint{\def\rm{\fam0\tenrm}
\textfont0=\tenrm \scriptfont0=\sevenrm \scriptscriptfont0=\fiverm
\textfont1=\teni  \scriptfont1=\seveni  \scriptscriptfont1=\fivei
\textfont2=\tensy \scriptfont2=\sevensy \scriptscriptfont2=\fivesy
\textfont\itfam=\tenit \def\it{\fam\itfam\tenit}\def\footnotefont{\ninepoint}%
\textfont\bffam=\tenbf \def\bf{\fam\bffam\tenbf}\def\sl{\fam\slfam\tensl}\rm}
\font\ninerm=cmr9 \font\sixrm=cmr6 \font\ninei=cmmi9 \font\sixi=cmmi6
\font\ninesy=cmsy9 \font\sixsy=cmsy6 \font\ninebf=cmbx9
\font\nineit=cmti9 \font\ninesl=cmsl9 \skewchar\ninei='177
\skewchar\sixi='177 \skewchar\ninesy='60 \skewchar\sixsy='60
\def\ninepoint{\def\rm{\fam0\ninerm}
\textfont0=\ninerm \scriptfont0=\sixrm \scriptscriptfont0=\fiverm
\textfont1=\ninei \scriptfont1=\sixi \scriptscriptfont1=\fivei
\textfont2=\ninesy \scriptfont2=\sixsy \scriptscriptfont2=\fivesy
\textfont\itfam=\ninei \def\it{\fam\itfam\nineit}\def\sl{\fam\slfam\ninesl}%
\textfont\bffam=\ninebf \def\bf{\fam\bffam\ninebf}\rm}
%
%

\hyphenation{anom-aly anom-alies coun-ter-term coun-ter-terms}
\def\inv{^{\raise.15ex\hbox{${\scriptscriptstyle -}$}\kern-.05em 1}}

\def\Dsl{\,\raise.15ex\hbox{/}\mkern-13.5mu D} 
\def\dsl{\raise.15ex\hbox{/}\kern-.57em\partial}

\font\bigit=cmti10 scaled \magstep1
\def\lspace{\ifx\answ\bigans{}\else\qquad\fi}
\def\lbspace{\ifx\answ\bigans{}\else\hskip-.2in\fi} 
\def\boxeqn#1{\vcenter{\vbox{\hrule\hbox{\vrule\kern3pt\vbox{\kern3pt
    \hbox{${\displaystyle #1}$}\kern3pt}\kern3pt\vrule}\hrule}}}
\def\mbox#1#2{\vcenter{\hrule \hbox{\vrule height#2in
        \kern#1in \vrule} \hrule}}  
%
 \def\CC{{\cal C}}

\def\darr#1{\raise1.5ex\hbox{$\leftrightarrow$}\mkern-16.5mu #1}

\def\roughly#1{\raise.3ex\hbox{$#1$\kern-.75em\lower1ex\hbox{$\sim$}}}

\def\frac#1#2{{#1\over#2}}

\def\journal#1&#2(#3){\unskip, #1~\bf #2 \rm(19#3) }
\def\andjournal#1&#2(#3){\sl #1~\bf #2 \rm (19#3) }

\def\bra#1{\left\langle #1\right|}
\def\ket#1{\left| #1\right\rangle}
\def\det{{\rm det}}

\catcode`\@=11\def\slash#1{\mathord{\mathpalette\c@ncel{#1}}}
\overfullrule=0pt
\def\steepslash{\c@ncel}
\def\frac#1#2{{#1\over #2}}

\def\:{\!:\!}
\def\inbar{\,\vrule height1.5ex width.4pt depth0pt}
\def\IQ{\relax\,\hbox{$\inbar\kern-.3em{\rm Q}$}}
\def\IB{\relax{\rm I\kern-.18em B}}
\def\IC{\relax\hbox{$\inbar\kern-.3em{\rm C}$}}
\def\IP{\relax{\rm I\kern-.18em P}}
\def\IR{\relax{\rm I\kern-.18em R}}
\def\ZZ{\relax\ifmmode\mathchoice
{\hbox{Z\kern-.4em Z}}{\hbox{Z\kern-.4em Z}}
{\lower.9pt\hbox{Z\kern-.4em Z}}
{\lower1.2pt\hbox{Z\kern-.4em Z}}\else{Z\kern-.4em Z}\fi}

\catcode`\@=12

\def\npb#1(#2)#3{{ Nucl. Phys. }{B#1} (#2) #3}
\def\plb#1(#2)#3{{ Phys. Lett. }{#1B} (#2) #3}
\def\pla#1(#2)#3{{ Phys. Lett. }{#1A} (#2) #3}
\def\prl#1(#2)#3{{ Phys. Rev. Lett. }{#1} (#2) #3}
\def\mpla#1(#2)#3{{ Mod. Phys. Lett. }{A#1} (#2) #3}
\def\ijmpa#1(#2)#3{{ Int. J. Mod. Phys. }{A#1} (#2) #3}
\def\cmp#1(#2)#3{{ Comm. Math. Phys. }{#1} (#2) #3}
\def\cqg#1(#2)#3{{ Class. Quantum Grav. }{#1} (#2) #3}
\def\jmp#1(#2)#3{{ J. Math. Phys. }{#1} (#2) #3}
\def\anp#1(#2)#3{{ Ann. Phys. }{#1} (#2) #3}
\def\prd#1(#2)#3{{ Phys. Rev. } {D{#1}} (#2) #3}
\def\ptp#1(#2)#3{{ Progr. Theor. Phys. }{#1} (#2) #3}
\def\aom#1(#2)#3{{ Ann. Math. }{#1} (#2) #3}

\def\bs{\bigskip}

\def\br{\buildrel}
\def\bra{\langle}
\def\ket{\rangle}

\def\C{{\bf C}}

\def\K{{\bf K}}

\def\P{{\bf P}}

\def\cA{{\cal A}}

\def\cE{{\cal E}}
\def\cF{{\cal F}}

\def\cH{{\cal H}}

\def\cK{{\cal K}}
\def\cL{{\cal L}}
\def\cM{{\cal M}}

\def\cO{{\cal O}}

\def\cR{{\cal R}}
\def\cS{{\cal S}}
\def\cT{{\cal T}}

\def\cV{{\cal V}}

\def\cX{{\cal X}}

\def\cicy#1(#2|#3)#4{\left(\matrix{#2}\right|\!\!
                     \left|\matrix{#3}\right)^{{#4}}_{#1}}

\def\lra{\longrightarrow}
\def\lla{\longleftarrow}
\def\ra{\rightarrow}

\def\bs{\bigskip}

\def\Box{{\,\lower0.9pt\vbox{\hrule
\hbox{\vrule height 0.2 cm \hskip 0.2 cm
\vrule height 0.2 cm}\hrule}\,}}

\global\newcount\thmno \global\thmno=0
\def\definition#1{\global\advance\thmno by1
\bigskip\noindent{\bf Definition \secsym\the\thmno. }{\it #1}
\par\nobreak\medskip\nobreak}
\def\question#1{\global\advance\thmno by1
\bigskip\noindent{\bf Question \secsym\the\thmno. }{\it #1}
\par\nobreak\medskip\nobreak}
\def\theorem#1{\global\advance\thmno by1
\bigskip\noindent{\bf Theorem \secsym\the\thmno. }{\it #1}
\par\nobreak\medskip\nobreak}
\def\proposition#1{\global\advance\thmno by1
\bigskip\noindent{\bf Proposition \secsym\the\thmno. }{\it #1}
\par\nobreak\medskip\nobreak}
\def\corollary#1{\global\advance\thmno by1
\bigskip\noindent{\bf Corollary \secsym\the\thmno. }{\it #1}
\par\nobreak\medskip\nobreak}
\def\lemma#1{\global\advance\thmno by1
\bigskip\noindent{\bf Lemma \secsym\the\thmno. }{\it #1}
\par\nobreak\medskip\nobreak}
\def\conjecture#1{\global\advance\thmno by1
\bigskip\noindent{\bf Conjecture \secsym\the\thmno. }{\it #1}
\par\nobreak\medskip\nobreak}
\def\exercise#1{\global\advance\thmno by1
\bigskip\noindent{\bf Exercise \secsym\the\thmno. }{\it #1}
\par\nobreak\medskip\nobreak}
\def\remark#1{\global\advance\thmno by1
\bigskip\noindent{\bf Remark \secsym\the\thmno. }{\it #1}
\par\nobreak\medskip\nobreak}
\def\problem#1{\global\advance\thmno by1
\bigskip\noindent{\bf Problem \secsym\the\thmno. }{\it #1}
\par\nobreak\medskip\nobreak}
\def\others#1#2{\global\advance\thmno by1
\bigskip\noindent{\bf #1 \secsym\the\thmno. }{\it #2}
\par\nobreak\medskip\nobreak}
\def\proof{\noindent Proof: }

\def\thmlab#1{\xdef #1{\secsym\the\thmno}\writedef{#1\leftbracket#1}\wrlabeL{#1=#1}}
%
%
\def\newsec#1{\global\advance\secno by1\message{(\the\secno. #1)}
\global\subsecno=0\thmno=0\eqnres@t\noindent{\bf\the\secno. #1}
\writetoca{{\secsym} {#1}}\par\nobreak\medskip\nobreak}
\def\eqnres@t{\xdef\secsym{\the\secno.}\global\meqno=1\bigbreak\bigskip}
\def\sequentialequations{\def\eqnres@t{\bigbreak}}\xdef\secsym{}
%

%


\ref\Bertram{A. Bertram, {\it Another way to enumerate rational curves
with torus action}, math.AG/9905159.}
\ref\BDPP{G. Bini, C. De Concini, M. Polito, and C. Procesi,
{\it Givental's work relative to mirror symmetry},
math.AG/9805097.}
\ref\Elezi{A. Elezi, {\it Mirror symmetry for concavex bundles
on projective spaces}, math.AG/0004157.}
\ref\FP{C. Faber, and R. Pandharipande, {\it Hodge Integrals and
Gromov-Witten Theory}, math.AG/9810173.}
\ref\FriedmanMorgan{R. Friedman, and J.W. Morgan,
{\it Smooth four-manifolds and complex surfaces,} in
Ergebnisse der Mathematik und ihrer Grenzgebiete (3)
[Results in Mathematics and Related Areas (3)], 27. Springer-Verlag, Berlin, 1994.}
\ref\Giv{A. Givental, {\it Equivariant Gromov-Witten
invariants}, alg-geom/9603021.}
\ref\GP{T. Graber and R. Pandharipande, {\it Localization of
virtual classes}, alg-geom/9708001.}
\ref\KKV{S. Katz, A. Klemm, and C. Vafa, {\it Geometric engineering
of quantum field theories}, Nucl. Phys. B497 (1997) 173-195.}
\ref\KM{F.F. Knudsen, and D. Mumford,
{\it The projectivity of the moduli space of stable curves.
I. Preliminaries on "det" and "Div". }
Math. Scand. 39 (1976), no. 1, 19-55.}
\ref\K{M. Kontsevich,
{\it Enumeration of rational curves via torus actions.}
In: The Moduli Space of Curves, ed. by
R. Dijkgraaf, C. Faber, G. van der Geer, Progress in Math.
vol. 129, Birkh\"auser, 1995, 335--368.}
\ref\LiTianII{ J. Li and G. Tian, {\it
Virtual moduli cycle and
Gromov-Witten invariants of algebraic varieties},
J. of Amer. math. Soc. 11, no. 1, (1998) 119-174.}
\ref\LLYI{B. Lian, K. Liu and S.T. Yau, {\it Mirror Principle I},
Asian J. Math. Vol. 1, No. 4 (1997) 729-763.}
\ref\LLYII{B. Lian, K. Liu and S.T. Yau, {\it Mirror Principle II},
Asian J. Math. Vol. 3, No. 1 (1999).}
\ref\LLYIII{B. Lian, K. Liu and S.T. Yau, {\it Mirror Principle III},
math.AG/9912038.}
\ref\Pandha{R. Pandharipande, {\it Rational curves on
hypersurfaces (after givental)},
math.AG/9806133.}

\Title{}{Mirror Principle IV}
 \centerline{\titlerms Bong H. Lian}
 \centerline{\it Department of Mathematics}                \vskip-1mm
 \centerline{\it Brandeis University, Waltham, MA 02154}   \vskip-1mm
 \centerline{ lian@brandeis.edu}         \vskip-1mm
\vskip .2in
\centerline{\titlerms Kefeng Liu}
 \centerline{\it Department of Mathematics}      \vskip-1mm
 \centerline{\it Stanford University, Stanford, CA 94305}  \vskip-1mm
 \centerline{ kefeng@math.stanford.edu}                   \vskip-1mm
\vskip .2in
 \centerline{\titlerms Shing-Tung Yau}
 \centerline{\it Department of Mathematics}                \vskip-1mm
 \centerline{\it Harvard University, Cambridge, MA 02138}  \vskip-1mm
 \centerline{ yau@math.harvard.edu}                  \vskip-1mm
\vfill

{\it In honor of Professors M. Atiyah, R. Bott, F. Hirzebruch, and I. Singer.}
\vskip .2in
Abstract.
This is a continuation  of {\it Mirror Principle III}
\LLYIII.


\Date{May 2000}

\newsec{Some Background}

This paper is a sequel to \LLYI\LLYII\LLYIII. In this series of papers
we develop {\it mirror principle} in increasing generality and breadth.
Given a projective manifold $X$, mirror principle is a theory
that yields relationships for and often computes the intersection
numbers of cohomology classes of the form
$b(V_D)$ on stable moduli spaces $\bar M_{g,k}(d,X)$.
Here $V_D$ is a certain induced vector bundles on $\bar M_{g,k}(d,X)$
and $b$ is any given multiplicative cohomology class.
In the first paper \LLYI, we consider this problem in the genus zero $g=0$ case
when $X=\P^n$ and $V_D$ is a bundle induced by any convex and/or concave bundle
$V$ on $\P^n$. As a consequence, we have proved a mirror formula which
computes the intersection numbers via a generating function.
When $X=\P^n$, $V$ is a direct sum of positive line bundles on $\P^n$,
and $b$ is the Euler class,
a second proof of this special case has been given
in \Pandha\BDPP~ following an approach proposed in \Giv.
Other proofs
in this case has also been given in \Bertram,
and when $V$ includes negative line bundles,
in \Elezi.
In \LLYII, we develop mirror principle
when $X$ is a projective manifold with $TX$ convex.
In \LLYIII, we consider the $g=0$ case when
$X$ is an arbitrary projective manifold.
Here  emphasis has been
put on a class of $T$-manifolds (which we call balloon manifolds)
because in this case mirror principle yields a (linear!)
reconstruction algorithm which
computes in principle all the intersection numbers above
for {\it any} convex/concave equivariant bundle $V$ on $X$ and {\it any}
equivariant multiplicative class $b$.
Moreover, specializing this theory to the case of
line bundles on toric manifolds and $b$ to Euler class,
we give a proof of the mirror formula for
toric manifolds.
In \LLYIII, we have also begun
to develop mirror principle for higher genus. In this paper,
we complete the theory for hight genus.
We also extend the
theory to include the intersection numbers for
cohomology classes of the form $ev^*(\phi)b(V_D)$.
Here $ev:\bar M_{g,k}(d,X)\ra X^k$ is the usual evaluation map
into the product $X^k$ of $k$ copies of $X$, and $\phi$ is any cohomology
class on $X^k$.

For motivations and the main ideas of mirror principle,
we refer the reader to the introduction of \LLYI\LLYII.

{\it Acknowledgement:} B.H.L. wishes to
thank the organizers for inviting him to lecture
in the Conference on Geometry
and Topology in honor of M. Atiyah, R. Bott, F. Hirzebruch,
and I. Singer. Once again, we owe our special thanks to
J. Li who has been patiently lending his help to us
throughout this project.
B.H.L.'s research
is supported by NSF grant DMS-0072158.
K.L.'s research is supported by NSF grant
DMS-9803234 and the Terman fellowship and the Sloane fellowship.
S.T.Y.'s research is supported by DOE grant
DE-FG02-88ER25065 and NSF grant DMS-9803347.

\newsec{Higher Genus}

We assume that the reader is familiar with \LLYIII.
We follow the notations introduced there. Most results proved
here have been summarized in section 5.5 there.

In the first subsection, we give some examples of gluing sequences,
a notion introduced in \LLYIII. We also prove a quadratic identity,
which is a generalization of Theorem 3.6 in \LLYIII~to higher genus
plus multiple marked points.
In the second subsection, we give a reconstruction algorithm
which allows us to reconstruct the Euler series $A(t)$ associated to
a gluing sequence in terms of Hodge integrals and some leading
terms of $A(t)$.

\item{1.} Throughout this note, we abbreviate
the data $(g,k;d)$ as $D$ and write
$$M_D=
M_{g,k}((d,1),X\times\P^1).$$
We denote by $LT_D(X),LT_{g,k}(d,X)$ the Li-Tian class of $M_D$
and $M_{g,k}(d,X)$ respectively.

\item{2.} In the last subsection, we prove the regularity of
the collapsing map
$$\varphi:M_{g,0}((d,1),\P^n\times \P^1)\ra N_d$$
generalizing Lemma 2.6 in \LLYI. Replacing $\P^n$ by
the product
$$Y=\P^{n_1}\times\cdots\times\P^{n_m},$$
and
$N_d$ by
$$W_d=N_{d_1}\times\cdots\times N_{d_m},$$
we get the map
$$M_{g,0}((d,1),Y\times \P^1)\ra W_d.$$
Given
an equivariant projective embedding $\tau:X\ra Y$ with $A^1(X)\cong A^1(Y)$
(see \LLYIII),
we have an induced map
$M_D=M_{g,k}((d,1),X\times\P^1)\ra M_{g,0}((d,1),Y\times\P^1)$.
Composing this with
$M_{g,0}((d,1),Y\times \P^1)\ra W_d$,
we get a $G$-equivariant map
$$M_D\ra W_d$$
which we also denoted by $\varphi$.
This map will be used substantially to do functorial localization in the
first subsection.

\item{3.} The standard $\C^\times$ action on $\P^1$
induces an action on each $M_D$ (see section 5.5 \LLYIII). The fixed point components are
labelled by
$F_{D_1,D_2}$ with $D_1=(g_1,k_1;d_1)$, $D_2=(g_2,k_2;d_2)$,
$g_1+g_2=g$, $k_1+k_2=k$, $d_1+d_2=d$.
As in the genus zero case, when $d_1,d_2\neq0$,
a stable map $(C,f,y_1,...,y_k)$ in this component is given by gluing
two pointed stable maps $(f_1,C_1,y_1,..,y_{k_1},x_1)\in M_{g_1,k_1+1}(d_1,X),~
(f_2,C_2,y_{k_1+1},...,y_k,x_2)\in M_{g_2,k_2+1}(d_2,X)$ with $f_1(x_1)=f_2(x_2)$,
to $C_0=\P^1$ at $0$ and $\infty$ at the marked points (cf. section 3).
We can therefore identify
$F_{D_1,D_2}$ with
$M_{g_1,k_1+1}(d_1,X)\times_X M_{g_2,k_2+1}(d_2,X)$.
We also have a special component $F_{D,O}$
which is obtained by gluing a $k+1$ pointed stable map to
$\P^1$ at either $0$, as described above. Likewise for $F_{O,D}$.
We denote by
$$
i:F_{D_1,D_2}\ra M_D,
$$
the inclusions.

\item{4.} There are two obvious projection maps
$$
p_0:F_{D_1,D_2}\ra M_{g_1,k_1+1}(d_1,X),~~
p_\infty:F_{D_1,D_2}\ra M_{g_2,k_2+1}(d_2,X).
$$
The map $p_0$ strips away (with the notations above)
the stable maps $(f_2,C_2,y_{k_1+1},...,y_k,x_2)$
glued to the $\P^1$ at $\infty$, and forgets the $\P^1$;
$p_\infty$ strips away the stable map
$(f_1,C_1,y_1,...,y_{k_1},x_1)$ glued to the $\P^1$ at $0$
and forgets the $\P^1$.

\item{5.}
We also have the following evaluation maps, and the forgetting map:
$$
e:F_{D_1,D_2}\ra X,~~~e_D:M_{g,k+1}(d,X)\ra X,~~~
\rho:M_{g,k+1}(d,X)\ra M_{g,k}(d,X).$$
Here $e$ evaluates a stable map in $F_{D_1,D_2}$ at the gluing point,
$e_D$ evaluates a $k+1$ pointed stable map at the last marked point,
and $\rho$ forgets the last marked point.
Relating and summarizing the natural maps above is the following
diagram:
\eqn\CommDiagram{
\matrix{
X & {\br e\over\lla} &  F_{D_1,D_2} &
{\br i\over\lra} &  M_D & {\br \pi\over\lra} & M_{g,k}(d,X) & \cr
e_{D_1}\uparrow & p_0\swarrow &   &\searrow p_\infty  & & & \cr
M_{g_1,k_1+1}(d_1,X) &   &   &    & M_{g_2,k_2+1}(d_2,X)  &  & \cr
\rho\downarrow &  &   &  &\downarrow\rho &  &  \cr
M_{g_1,k_1}(d_1,X) &   &   &    & M_{g_2,k_2}(d_2,X)  &  &
}
}
Here $\pi$ is natural morphism which maps $(C,f,y_1,..,y_k)\in M_D$
to the stabilization of $(C,\pi_1\circ f,y_1,..,y_k)$, where $\pi_1:X\times\P^1\ra X$
is the projection. Note that
we can identify $M_{g,k+1}(d,X)$ with $F_{O,D}$ via $p_\infty$.
When $D_1=O:=(0,0;0)$,
the right part of the diagram above completes to a commutative triangle, ie.
$\pi\circ i=\rho$.

\item{6.}
Let $\bar M_{g,k}$ be the Deligne-Mumford moduli space
of $k$-pointed, genus $g$ stable curves.
Recall the map
$$M_{g,k}(d,X)\ra \bar M_{g,k}$$
which sends $(C,f,y_1,..,y_k)$ to the stabilization of $(C,y_1,...,y_k)$.
Let $\cL$ and $\cH$ be respectively the universal line bundle and
the Hodge bundle on $\bar M_{g,k+1}$. Thus $\cL,\cH$ have
fibers at $(C,y_1,...,y_k,x)$ given by $T_xC$ and
$H^0(C,K_C)\cong H^1(C,\cO)^*$ respectively.
We denote
$$\Lambda_g(\xi)=\sum_{i=0}^g\xi^i c_{g-i}(\cH)$$
for any formal variable $\xi$ (cf. \FP). We denote by the same
notation the pullback of $\Lambda_g(\xi)$ to $M_{g,k+1}(d,X)$.
We denote by $L_D$ the universal line bundle on $M_{g,k+1}(d,X)$.
corresponding to the last marked point.

\item{7.}
{\it Notations.} In all formulas below involving $d,d_1,d_2$, and $g,g_1,g_2$,
it is always assumed that
$$\eqalign{
D&=(g,k;d)\cr
D_1&=(g_1,k_1;d_1)\cr
D_2&=(g_2,k_2;d_2)\cr
g&=g_1+g_2\cr
k&=k_1+k_2\cr
d&=d_1+d_2.
}
$$

\lemma{(cf. Lemma 3.5 \LLYIII)
Let $g=g_1+g_2$, $k_1+k_2=k$, $d=d_1+d_2$.
For $d_1,d_2\neq 0$, the equivariant Euler class of
the virtual normal bundle $N_{F_{D_1,D_2}/M_D}$ is
$$
e_G(F_{D_1,D_2}/M_D)=
p_0^*\left(\alpha(\alpha+c_1(L_{D_1}))~\Lambda_{g_1}(\alpha)^{-1}\right)
p_\infty^*\left(
\alpha(\alpha-c_1(L_{D_2}))~\Lambda_{g_2}(-\alpha)^{-1}\right).
$$
For $d_1=0,~g_1=0,~k_1=0$,
$$
e_G(F_{O,D}/M_D)
=-\alpha(-\alpha+c_1(L_D))~\Lambda_g(-\alpha)^{-1}.
$$
For $d_2=0,~g_2=0,~k_2=0$,
$$
e_G(F_{D,O}/M_D)
=\alpha(\alpha+c_1(L_D))~\Lambda_g(\alpha)^{-1}.
$$
}
\thmlab\NormalBundle
\proof
We consider the first equality,
the other two being similar.
We will compute the virtual normal bundle $N_{F_{D_1,D_2}/M_D}$
of $F_{D_1,D_2}$ in $M_D$ following the methods in \K\GP\LLYI,
using the description of $F_{D_1,D_2}$ given above.
We must identify the terms
appearing in the tangent obstruction sequence of $M_D$. (See
\GP~section 4.)

Consider the
bundle $V:=\pi_1^*TX\oplus \pi_2^* T\P^1$ on $X\times\P^1$, where $\pi_1$ and
$\pi_2$ are the projections from $X\times \P^1$ to $X$ and
$\P^1$ respectively.

According to the description above,
for each stable map
$(C,f,x_1,...,x_k)$ in $F_{D_1,D_2}$,
we have an exact sequence over $C$:
\eqn\dumb{
0\rightarrow f^*V\rightarrow f_1^*V\oplus f_0^*V\oplus f_2^*V
\rightarrow V_{x_1}\oplus
V_{x_2}\rightarrow 0.
}
Here $f_0$ is the restriction of $f$ to $C_0$, and
$V_{x_1}, \ V_{x_2}$ denote respectively the
bundles
$$
\pi_1^*T_xX \oplus \pi_2^* T_0\P^1,
\ \ \pi_1^*T_xX \oplus \pi_2^*
T_\infty\P^1
$$
where $x=f_0(C_0)\in X$.
From the  long exact sequence associated to \dumb,
we get the equality in the K-group:
$$\eqalign{
&H^0(C, f^*V)-H^1(C, f^*V)\cr
&~~~=\sum_{j=0}^2 [H^0(C_j, f_j^*V)-H^1(C_j, f_j^*V)]
-(V_{x_1}+V_{x_2}).
}$$

The tangent complex of $M_D$ restricted to
$F_{D_1,D_2}$ is
\eqn\Total{\eqalign{
&H^0(C, f^*V)-H^1(C, f^*V)
+T_{x_1}C_1\otimes T_0\P^1+
T_{x_2}C_2\otimes T_\infty\P^1-A_C\cr
&=\sum_{j=0}^2 [H^0(C_j, f_j^*V)-H^1(C_j, f_j^*V)]
-(V_{x_1}+V_{x_2})\cr
&+T_{x_1}C_1\otimes T_0\P^1+
T_{x_2}C_2\otimes T_\infty\P^1-A_C.
}}
where
$T_{x_1}C_1\otimes T_0\P^1$ and
$T_{x_2}C_2\otimes T_\infty\P^1$ are contributions
from the deformation of the nodes at $x_1,x_2$ of $C$, and
$A_C$ is the contribution from the infinitesimal automorphisms
of $C$.
To get the moving parts of this,
we subtract from it the fixed parts corresponding to
$F_{D_1,D_2}=M_{g_1,k_1+1}(d_1,X)\times_X
M_{g_2,k_2+1}(d_2,X)$ (see description above). This is
given by
\eqn\Fixed{
\sum_{j=1}^2[H^0(C_j,f^*_j\pi_1^*TX)-H^1(C_j,f^*_j\pi_1^*TX)]-\pi^*_1T_xX-A_{C'},
}
where $C'$ is the curve obtained from $C$ by contracting the component $C_0$.
Note that here we have ignored the contributation coming from
the deformation of $C_1,C_2$ in both \Total~ and \Fixed,
because the same contributation appear in both and
hence this contributation cancels out in the difference.

We now compute and compare the terms in both \Total~ and \Fixed~ above.
Since $f_0$ maps $C_0$ to a point $x$, we have
$$\eqalign{
H^1(C_0, f_0^*\pi_1^*TX)&\simeq H^1(C_0, O)\otimes \pi_1^*TX =0\cr
H^1(C_0, f_0^*\pi_2^*T\P^1)&\simeq H^1(C_0, TC_0) =0\cr
H^0(C_0, f_0^*\pi_1^*TX)&=\pi_1^*T_xX.
}$$
Similarly for $C_1$ and $f_1$, we have
$$\eqalign{
H^1(C_1, f_1^*V)&=H^1(C_1, f_1^*\pi_1^*TX)+H^1(C_1, f_1^*\pi_2^*T\P^1)\cr
H^1(C_1, f_1^*\pi_2^*T\P^1)&\simeq H^1(C_1, O)\otimes \pi_2^* T_0\P^1\cr
H^0(C_1, f_1^*V)&=H^0(C_1, f_1^*\pi_1^*TX)+H^0(C_1, f_1^*\pi_2^*T\P^1)\cr
H^0(C_1, f_1^*\pi_2^*T\P^1)&\simeq \pi_2^*T_0\P^1.
}$$
Likewise, we have similar relations with $C_1,f_1, T_0\P^1$
replaced by $C_2,f_2,T_\infty\P^1$ everywhere.


Putting these formulas together, we get
$$\eqalign{
N_{F_{D_1,D_2}/M_D}
&= H^0(C_0, f_0^*\pi_2^*T\P^1)+T_{x_1}C_1\otimes T_0\P^1\cr
&~~+T_{x_2}C_2\otimes T_\infty\P^1 -H^1(C_1, O)\otimes
\pi_2^*T_0\P^1-H^1(C_2, O)\otimes \pi_2^*T_\infty\P^1-A_{C_0}.
}$$
By taking equivariant Euler classes, we get the required formula.
Here the terms
$-H^1(C_1, O)\otimes
\pi_2^*T_0\P^1$ and $-H^1(C_2, O)\otimes \pi_2^*T_\infty\P^1$
contribute $p_0^*\Lambda_{g_1}(\alpha)^{-1}$ and
$p^*_\infty\Lambda_{g_2}(-\alpha)^{-1}$
respectively. The terms
$T_{x_1}C_1\otimes T_0\P^1$ and
$T_{x_2}C_2\otimes T_\infty\P^1$ contribute
$\alpha+p_0^*c_1(L_{D_1})$ and $-\alpha+p_\infty^*c_1(L_{D_2})$
respectively, and the term
$H^0(C_0, f_0^*\pi_2^*T\P^1)-A_{C_0}$ contributes $-\alpha^2$
(see \LLYI, section 2.3).

Similarly, when $d_1=0, g_1=0, k_1=0$, we have
$$
N_{F_{O,D}/M_D} = H^0(C_0, f_0^*\pi_2^*T\P^1)+T_{x_2}C_2\otimes
T_\infty\P^1-H^1(C_2, O)\otimes \pi_2^*T_\infty\P^1-A_{C_0}.
$$
The term $-A_{C_0}$ contributes a factor $-\alpha$.
Similarly for $d_2=0$ and $g_2=0$, we have
$$
N_{F_{D,O}/M_D} = H^0(C_0, f_0^*\pi_2^*T\P^1)+T_{x_1}C_1\otimes
T_0\P^1-H^1(C_1, O)\otimes \pi_2^*T_0\P^1-A_{C_0},
$$
and now the term $-A_{C_0}$
contributes a factor $\alpha$. $\Box$

\subsec{Gluing sequences}

Fix a class $\Omega\in A_T^*(X)$, such that $\Omega^{-1}$ is well-defined.
We call the list of classes
$$
b_D\in A_T^*(M_{g,k}(d,X)).
$$
an $\Omega$-gluing sequence
if we have the (gluing) identities on the $F_{D_1,D_2}$:
$$
e^*\Omega\cdot i^*\pi^* b_D
=p^*_0\rho^* b_{D_1}\cdot
p^*_\infty\rho^* b_{D_2}.
$$
(This generalizes the definition in \LLYIII~ to include
the cases with multiple marked points.)
Gluing sequences have the following
obvious multiplicative property: that if $b_D$ and $b'_D$ form
two gluing sequences with respective to, say $\Omega$ and $\Omega'$,
then the product $b_D b'_D$ form a gluing sequence
with respective to $\Omega\Omega'$.

Let $V$ be a $T$-equivariant bundle on $X$. Suppose that
$H^0(C,f^*V)=0$ for every positive degree map $f:C\ra X$ where $C$
a nonsingular genus $g$ curve. Then $V$ induces on each $M_{g,k}(d,X)$
a vector bundle $V_D$ whose fiber at $(C,f,y_1,..,y_k)$ is $H^1(C,f^*V)$.
We call such a $V$ a {\it concave bundle} on $X$.

\bs
\noindent {\it Example 1.} $X=\P^n$, and $V=\cO(-k)$, $k<0$.
\bs

\noindent {\it Example 2.} If $X$ is a projective manifold with
$V=K_X<0$, then $V$ induces the bundles $V_D$. This is
the situation in {\it local mirror symmetry} \LLYI\KKV.
\bs

Let $b_T$ be a $T$-equivariant multiplicative class.
such that $\Omega=b_T(V)^{-1}$ is well-defined.
Consider the classes $b_T(V_D)\in A_T^*(M_{g,k}(d,X))$.

\lemma{The cohomology classes $b_T(V_D)$
form an $\Omega$-gluing sequence.}
\proof
The proof is essentially the same as the argument for the genus
zero gluing identity for a concave bundle $V$. See
the first half of the proof of Theorem 3.6 \LLYIII. $\Box$

We now discuss a second example of a gluing sequence.
Recall that a point $(f,C)$ in
$F_{D_1,D_2}$ comes from gluing together a pair of
stable maps
$(f_1,C_1,y_1,..,y_{k_1},x_1),(f_2,C_2,y_{k_1+1},..,y_k,x_2)$ with $f_1(x_1)=f_2(x_2)$,
to $C_0=\P^1$ at $0$ and $\infty$, so that
we have a long exact sequence
$$\eqalign{
&0\ra H^0(C,\cO)\ra H^0(C_1,\cO)\oplus H^0(C_2,\cO)\ra H^0(C_0,\cO)\cr
&\ra H^1(C,\cO)\ra H^1(C_1,\cO)\oplus H^1(C_2,\cO)\ra H^1(C_0,\cO)\ra 0.
}$$
Thus we have a natural isomorphism
$$
H^1(C,\cO)\cong H^1(C_1,\cO)\oplus H^1(C_2,\cO).
$$
This implies the isomorphism
$$i^*\pi^*\cH_D=p_0^* \rho^*\cH_{D_1}\oplus p_\infty^*\rho^*\cH_{D_2}$$
of bundles on $F_{D_1,D_2}$.
Here $\cH_D$ denotes the bundle on $M_{g,k}(d,X)$
with fiber $H^1(C,\cO)$. (Note that for $g\geq2$,
then $\cH_D$ is the pullback of the dual of the Hodge bundle $\cH$ via
the natural map $M_{g,k}(d,X)\ra \bar M_{g,k}$.)
Thus if $b$ is a multiplicative class, and $b_D:=b(\cH_D)$, then
the isomorphism above yields the gluing identity
$$
i^*\pi^*b_D= p_0^* \rho^*b_{D_1}\cdot
p_\infty^*\rho^*b_{D_2}
$$
with $\Omega=1$. To summarize, we have

\lemma{The cohomology classes $b_D:=b(\cH_D)$
above form a $1$-gluing sequence.}

Note that in both examples above, each class $b_D$
is naturally the pullback, via the forgetful map,
of a class $b_d^g\in A_T^*(M_{g,0}(d,X))$. We will call
{\it a list of classes $b_d^g\in A_T^*(M_{g,0}(d,X))$ an $\Omega$-gluing
sequence, if their pullbacks to $M_{g,k}(d,X)$ form
a gluing sequence in the sense introduced above.}

We now discuss a third construction.
Fix a set of generators $\phi_i$ of $A_T^*(X)$,
as a free module over $A_T^*(pt)$. Let
$$\phi=\sum_i s_i\phi_i$$
where $s_i$ are formal variables. Let
$$
\phi_D:=\prod_{j=1}^k ev_j^*\phi\in A_T^*(M_{g,k}(d,X)).
$$
Here the map $ev_j$ evaluates at the $j$th marked point
of a stable map in $M_{g,k}(d,X)$. Then we have
$$
ev_j\circ\pi\circ i=ev_j\circ\rho\circ p_0
$$
for $j=1,..,k_1$. Here $ev_j$ on the LHS are evaluations on $M_{g,k}(d,X)$,
while $ev_j$ on the RHS are evaluations on $M_{g_1,k_1}(d_1,X)$.
Likewise we have
$$
ev_{k_1+j}\circ\pi\circ i=ev_j\circ\rho\circ p_\infty
$$
for $j=1,..,k_2$. It follows that
$$
i^*\pi^*\phi_D= p_0^* \rho^*\phi_{D_1}\cdot
p_\infty^*\rho^*\phi_{D_2}.
$$
Thus we have
\lemma{The cohomology classes $\phi_D$ above form a $1$-gluing sequence.}

Combining with the multiplicative property of gluing sequences,
as explained above,
this construction allows us to consider
the intersection numbers
of classes of the form $ev^*(\phi)~b(V_D)$ on stable map moduli.
In particular, this yields the
GW-invariants twisted by a multiplicative
class of the form $b(V_D)$.
Here $ev$ is the evaluation map $M_{g,k}(d,X)\ra X^k$
into the product of $k$ copies of $X$.
The results below are easily
generalized to the cases involving the additional factor $ev^*(\phi)$.

For $\omega\in A_G^*(M_D)$, introduce the notation (cf. section 3.2 \LLYIII)
$$
J_{D_1,D_2}\omega:=e_*
\left({i^*\omega\cap [F_{D_1,D_2}]^{vir}
\over e_G(F_{D_1,D_2}/M_D)}\right)\in A_*^T(X)(\alpha).
$$

\theorem{(cf. Theorem 3.6 \LLYIII)
Given a gluing sequence $b_D\in A_T^*(M_{g,k}(d,X))$, we have the
following identities in $A^T_*(X)(\alpha)$:
$$
\Omega\cap
J_{D_1,D_2}\pi^*b_D
=\overline{J_{O,D_1}\pi^*b_{D_1} }
\cdot J_{O,D_2}\pi^*b_{D_2}.
$$
}
\thmlab\GluingTheorem
\proof
Consider the fiber square
\eqn\FiberSquare{\matrix{
F_{D_1,D_2} & {\br \Delta'\over\longrightarrow}&
M_{g_1,k_1+1}(d_1,X)\times M_{g_2,k_2+1}(d_2,X)\cr
e\downarrow &  &~~~~~~~\downarrow e_{D_1}\times e_{D_2}\cr
X & {\br \Delta\over\longrightarrow} & X\times X
}}
where $\Delta$ is the diagonal inclusion.
Recall also that (section 6 \LiTianII)
$$
[F_{D_1,D_2}]^{vir}=\Delta^!
(LT_{g_1,k_1+1}(d_1,X)\times LT_{g_2,k_2+1}(d_2,X)).
$$
Put
$$
\omega={\rho^*b_{D_1}\over e_G(F_{D_1,O}/M_{D_1})}
\times{\rho^*b_{D_2}\over e_G(F_{O,D_2}/M_{D_2})}\cap
LT_{g_1,k_1+1}(d_1,X)\times LT_{g_2,k_2+1}(d_2,X).
$$
From the fiber square \FiberSquare, we have
$$
e_*\Delta^!(\omega)=\Delta^*(e_{D_1}\times e_{D_2})_*(\omega).
$$
On the one hand is
$$\eqalign{
\Delta^*(e_{D_1}\times e_{D_2})_*(\omega)
&={e_{D_1}}_*
{\rho^*b_{D_1}\cap LT_{g_1,k_1+1}(d_1,X)
\over e_G(F_{D_1,O}/M_{D_1})}\cdot
{e_{D_2}}_*
{\rho^*b_{D_2}\cap LT_{g_2,k_2+1}(d_2,X)
\over e_G(F_{O,D_2}/M_{D_2})}\cr
&=\overline{J_{O,D_1}\pi^*b_{D_1} }
\cdot J_{O,D_2}\pi^*b_{D_2}.
}$$
Here we have use the fact that $\pi\circ i=\rho$.
On the other hand, applying the gluing identity and Lemma \NormalBundle, we have
$$\eqalign{
e_*\Delta^!(\omega)&=
e_*\left(
p_0^*{\rho^* b_{D_1}\over e_G(F_{D_1,O}/M_{D_1})}
\cdot
p_\infty^*{\rho^* b_{D_2}\over e_G(F_{O,D_2}/M_{D_2})}
\cap [F_{D_1,D_2}]^{vir}\right)\cr
&=e_*\left({e^*\Omega\cdot i^*\pi^*b_D
\cap [F_{D_1,D_2}]^{vir}\over
e_G(F_{D_1,D_2}/M_D)}
\right)\cr
&=\Omega\cap J_{D_1,D_2}\pi^* b_D.
}
$$
This proves our assertion. $\Box$

\lemma{(cf. Lemma 3.2 \LLYIII)
Given a cohomology class $\omega$ on $M_D$, we have
the following identities on
the $\C^\times$ fixed point component $Y_{d_1,d_2}\cong Y$
in $W_d$:
$$
{j^*\varphi_*(\omega\cap LT_D(X))\over e_G(Y_{d_1,d_2}/W_d)}
=\sum_{g_1+g_2=g,k_1+k_2=k} \tau_*e_*\left({i^*\omega\cap[F_{D_1,D_2}]^{vir}
\over e_G(F_{D_1,D_2}/M_D)}\right).
$$
}
\thmlab\SquareLemma
\proof
This follows from applying functorial localization to the diagram
\eqn\CommSquareII{
\matrix{
\{F_{D_1,D_2}\} & {\br i\over \longrightarrow} & M_D\cr
\tau\circ e\downarrow &  & \downarrow \varphi\cr
Y_{d_1,d_2} & {\br j\over \longrightarrow} & W_d.}~~~~\Box
}

Now given a gluing sequence $b_D$, we put
$$
A_D:=J_{O,D}\pi^*b_D,~~~ A_d:=\sum_{g,k} A_D~ \nu^g \mu^k,~~~
A(t):=e^{-H\cdot t/\alpha}\sum_d A_d~e^{d\cdot t}.
$$
Here $\nu,\mu$ are formal variables. Consider the class
$\beta=\varphi_*(\pi^* b_D\cap LT_D(X))$. We have
$$\eqalign{
\int_{W_d} \beta\cap e^{\kappa\cdot\zeta}
&=\sum_{d_1+d_2=d}\int_{Y_{d_1,d_2}}
{j^*\beta\over e_G(Y_{d_1,d_2}/W_d)}~
e^{(H+d_1\alpha)\cdot\zeta}\cr
&=\sum_{D_1+D_2=D}\int_{Y_{d_1,d_2}}
\tau_* J_{D_1,D_2}\pi^*b_D~
e^{(H+d_1\alpha)\cdot\zeta}
~~~~~~~(Lemma~\SquareLemma)\cr
&=\sum_{D_1+D_2=D}\int_X
J_{D_1,D_2}\pi^*b_D~
e^{(H+d_1\alpha)\cdot\zeta}\cr
&=\sum_{D_1+D_2=D}\int_X
\Omega^{-1}\cap\overline{A_{D_1}} \cdot A_{D_2}~
e^{(H+d_1\alpha)\cdot\zeta}
~~~~~~(Theorem~\GluingTheorem).
}$$
Since $\beta\in A^G_*(W_d)$, hence
$\int_{W_d} \beta\cap c\in A_*^G(pt)=\C[\cT^*,\alpha]$
for all $c\in A_G^*(W_d)$, it follows that both sides of
the eqn. above lie in $\cR[[\zeta]]$. Thus we get

\corollary{$A(t)$
is an Euler series.}

We call $V$ be a $D$-critical concave bundle if
the homogeneous degree of the class $b(V_D)$
is the same as the expected dimension
of $M_{g,k}(d,X)$. We denote
$$
K_D=\int_{LT_{g,k}(d,X)} b(V_D).
$$

\lemma{If $V$ is a $D$-critical concave bundle, then
in the $T$-nonequivariant limit we have the following formula
$$
\int_X e^{-H\cdot t/\alpha} J_{O,D}\pi^* b(V_D)=(-1)^g
\alpha^{g-3}(2-2g-k-d\cdot t)K_D.
$$
}
\proof
The integral above is equal to
$$\eqalign{
&\int_{LT_{g, k+1}(d, X)}e^{-e^*H\cdot t/\alpha}
\frac{\rho^*b(V_D)\Lambda_g(-\alpha)}{\alpha(\alpha-c)}\cr
&=\int_{LT_{g, k}(d, X)}b(V_D)\Lambda_g(-\alpha)\rho_*
\left(\frac{e^{-e^*H\cdot t/\alpha}}
{\alpha(\alpha-c)}\right).
}$$
Here $\rho:M_{g,k+1}(d,X)\ra M_{g,k}(d,X)$ forgets the last marked point.

Since the fiber of $\rho$ is of dimension $1$, we take the degree $1$ term in the fiber
integral. The integration along the fiber $E$ is done
in essentially the same  way as in
the genus zero case (see Theorem 3.12(ii) of \LLYIII). It yields
$$
\int_E e^*H=d
$$
and
$$\int_E c =2-2g-k
$$
by Gauss-Bonnet formula. Since the degree of $b(V_D)$ coincides
with the dimension of $LT_{g,k}(d,X)$,
it follows that only the $\alpha^0$ term of $\Lambda_g(-\alpha)$
contributes to the integral above.
$\Box$

\subsec{Reconstruction}

From now on, we assume that $X$ is a balloon manifold,
as in sections 4.1 and 5.3 of \LLYIII.
We will find further constraints
on a gluing sequence by computing the linking
values of the Euler series $A(t)$.
Recall that in genus zero, the linking values of an Euler series,
say coming from $b_T(V_d)$, are determined by the
restrictions $i_F^*b_T(\rho^* V_d)$ to the
isolated fixed point $F=(\P^1,f_\delta,0)\in M_{0,1}(d,X)$,
which is a $\delta$-fold cover of a balloon $pq$ in $X$.
In higher genus with multiple marked points, this will be replaced by
a component in $M_{g,k+1}(d,X)$ consisting of
the following stable maps $(C,f,y_1,..,y_k,x)$.
Here $C$ is a union of two curves $C_1$ and
$C_0\cong\P^1$ such that $y_1,..,y_k\in C_1$ and
that $C_0{\br f\over\ra} pq$ is a $\delta$-fold
cover with $f(x)=p$, $f(C_1)=q$.
Therefore this fixed point component can be
identified canonically with
$\bar M_{g,k+1}$. For clarity, we will restrict the following
discussion to the case of $k=0$.
We'll denote the component by $F^g$. By convention,
$F^0$ is the isolated fixed point $(\P^1,f_\delta,0)$.
Recall that (section 5.3 \LLYIII)
\eqn\GenusZero{
e_G(F^0/M_{0,1}(d,X))^{-1}={-\lambda\over\delta}
{e_T[H^1(\P^1,f_\delta^*TX)]'\over
e_T[H^0(\P^1,f_\delta^*TX)]'}.
}

\theorem{Suppose $g>0$. Let $p\in X^T$,
$\omega\in A^*_T(M_{g,1}(d,X))[\alpha]$, and consider
$i_p^* e_*\left(
{\omega\cap LT_{g,1}(d,X)\over e_G(F_{O,D}/M_D)}\right)\in\C(\cT^*)(\alpha)$
as a function of $\alpha$. Then
\item{(i)} Every possible pole of the function is a scalar multiple of
a weight on $T_pX$.
\item{(ii)} Let $pq$ be a balloon in $X$,
and $\lambda$ be the weight on the tangent line $T_p(pq)$.
If $d=\delta[pq]\succ0$,
then the pole of the function at $\alpha=\lambda/\delta$
is of the form
$$
e_T(p/X)~{1\over\alpha(\alpha-\lambda/\delta)}
{1\over e_T(F^0/M_{0,1}(d,X))}
\int_{[F^g]^{vir}}{i_{F^g}^*\omega~\Lambda_g(\alpha)~ e_T(\cH^*\otimes T_qX)\over
(-{\lambda\over\delta}+c_1(\cL))}.
$$
}
\thmlab\SimplePole
\proof
The proof here is a slight modification of the genus zero case.
We repeat the details here for the readers' convenience.
Consider the commutative diagram
$$
\matrix{
\{F\} & {\br i_F\over\longrightarrow}&
M_{g,1}(d,X)\cr
e'\downarrow &  &~~e\downarrow\cr
p & {\br i_p\over\longrightarrow} & X
}
$$
where $e$ is the evaluation map,
$\{F\}$ are the fixed point components in $e^{-1}(p)$,
$e'$ is the restriction of $e$ to $\{F\}$, and $i_F,i_p$
are the usual inclusions. By functorial localization
we have, for any $\beta\in A^*_T(M_{g,1}(d,X))(\alpha)$,
\eqn\dumb{\eqalign{
i_p^* e_*(\beta\cap LT_{g,1}(d,X))
&=e_T(p/X)
~\sum_F e'_*\left(
{i_F^*\beta\cap[F]^{vir}\over e_T(F/M_{g,1}(d,X))}\right)\cr
&=e_T(p/X)~
\sum_F\int_{[F]^{vir}}
{i_F^*\beta\over e_T(F/M_{g,1}(d,X))}.
}}
We apply this to the class
$$
\beta=
{\omega\over e_G(F_{O,D}/M_D)}
={\omega~\Lambda_g(-\alpha)\over \alpha(\alpha-c)}
$$
where $c=c_1(L_D)$ (cf. Lemma \NormalBundle).
For (i), we will show that a pole of the sum \dumb~
is at either $\alpha=0$ or $\alpha=\lambda'/\delta'$
for some tangent weight $\lambda'$ on $T_pX$.
For (ii), we will show that only one component $F$ in the sum \dumb~
contributes to the pole at $\alpha=\lambda/\delta$,
that the contributing component is $F^g$,
and that the contribution has
the desired form.

A fixed point $(C,f,x)$ in $e^{-1}(p)$ is such that
$f(x)=p$, and that the image curve $f(C)$
lies in the union of the $T$-invariant
balloons in $X$.
The restriction of the first Chern class $c$ to
an $F$ must be of the form
$$i_F^*c=c_F+w_F$$
where $c_F\in A^1(F)$, and $w_F\in\cT^*$
is the weight of the representation on the line $T_xC$
induced by the linear map $df_x:T_xC\ra T_pX$.
The image of $df_x$ is either 0 or a tangent line $T_p(pr)$
of a balloon $pr$.
Thus $w_F$ is either zero (when the branch $C_1\subset C$ containing $x$
is contracted), or $w_F=\lambda'/\delta'$
(when $C_1{\br f\over\ra} X$ maps by a $\delta'$-fold cover
of a balloon $pr$ with tangent weight $\lambda'$).
The class $e_T(F/M_{g,1}(d,X))$
is obviously independent of $\alpha$.
Since $c_F$ is nilpotent, a pole of the sum \dumb~
is either at $\alpha=0$ or $\alpha=w_F$ for some $F$.
This proves (i).

Now, an $F$ in the
sum \dumb~ contributes to the pole at $\alpha=\lambda/\delta$
only if $w_F=\lambda/\delta$.
Since the weights on $T_pX$ are pairwise linearly independent,
that $\lambda/\delta=\lambda'/\delta'$ implies that
$\lambda=\lambda'$ and $\delta=\delta'$.
Since $d=\delta[pq]$, a fixed point $(C,f,x)$ contributing to
the pole at $\alpha=\lambda/\delta$ must have the following form:
that there is a branch $C_0\cong\P^1$ in $C$
such that $f|_{C_0}:C_0\ra pq$
is a $\delta$-fold
cover with $f(x)=p$. Let $y\in C_0$ be the preimage of $q$
under this covering.
Then the curve $C$ is a union of $C_0$
and a genus $g$ curve $C_1$ meeting $C_0$ at $y$,
and $f(C_1)=q$. In other words, the fixed point component $F$
contributing to
the pole at $\alpha=\lambda/\delta$ is $F^g\cong \bar M_{g,1}$.
It contributes to the sum \dumb~ the term
$$
\int_{[F^g]^{vir}}{i_{F^g}^*\beta\over e_T(F^g/M_{g,1}(d,X))}
=\int_{F^g} {i_{F^g}^*\omega~\Lambda_g(-\alpha)\over
\alpha(\alpha-\lambda/\delta-c_{F^g})} {1\over
e_T(F^g/M_{g,1}(d,X))}.
$$
Here $c_{F^g}\in A^1(F^g)$ is zero because the
universal line bundle $L^g_d$ restricted to $F^g$ is trivial
(the line $T_xC$ is located at the marked point $x$).

We now compute the virtual normal bundle $N_{F^g/M_{g,1}(d,X)}$.
A point $(C,f,x)$ in $F^g$ can be viewed as
gluing two stable maps $(C_0,f_0,x,y)\in M_{0,2}(d,X)$,
$(C_1,f_1,x_1)\in M_{g,1}(0,X)$,
by identifying $x_1\equiv y$.
Here $f_0:C_0\ra pq$
is a $\delta$-fold cover with $f_0(x)=p$, $f_0(y)=q$,
and $f_1(C_1)=q$.
As before, we have
$$\eqalign{
N_{F^g/M_{g,1}(d,X)}
&=[H^0(C,f^*TX)]-[H^1(C,f^*TX)]+[T_yC_1\otimes T_y C_0]
-A_{C_0}\cr
&=\left([H^0(C_0,f_0^*TX)]-H^1(C_0,f_0^*TX)]-A_{C_0}\right)\cr
&~~
-[H^1(C_1,f_1^*TX)]+[T_yC_1\otimes T_y C_0].
}$$
Note that the first three terms collected in the parentheses
is the virtual normal bundle of $F^0$ in $M_{0,1}(d,X)$.
The Euler class of this is constant on $F^g$, as
given in eqn. \GenusZero.
Since $f_1:C_1\ra X$ maps to the point $q$, it follows that
$[H^1(C_1,f_1^*TX)]=\cH^*\otimes T_qX$ where $\cH$ is the Hodge bundle
on $\bar M_{g,1}$. Clearly
$[T_yC_1\otimes T_y C_0]=\cL\otimes[{-\lambda\over\delta}]$,
where $\cL$ is the universal line bundle on $\bar M_{g,1}$,
and $[{-\lambda\over\delta}]$ is a 1 dimensional representation
of that given weight.
Therefore, we get
$$
e_T(F^g/M_{g,1}(d,X))=e_T(F^0/M_{0,1}(d,X))~
(-{\lambda\over\delta}+c_1(\cL))
~e_T(\cH^*\otimes T_qX)^{-1}.
$$
Hence the contribution of the sum \dumb~ to the pole
at $\alpha=\lambda/\delta$ is
$$\eqalign{
&e_T(p/X)~\int_{[F^g]^{vir}}{i_{F^g}^*\beta\over e_T(F^g/M_{g,1}(d,X))}\cr
&~~~~=e_T(p/X)~{1\over\alpha(\alpha-\lambda/\delta)}
{1\over e_T(F^0/M_{0,1}(d,X))}
\int_{[F^g]^{vir}}{i_{F^g}^*\omega~\Lambda_g(-\alpha) ~e_T(\cH^*\otimes T_qX)\over
{-\lambda\over\delta}+c_1(\cL))}.
}$$
This proves (ii).  $\Box$

Let $V$ be a concave bundle on $X$, and $b_T$ a choice of
multiplicative class as before.
Define the genus $g\geq0$,
degree $d=\delta[pq]$,
{\it linking values} at the balloon $pq$:
$$
Lk_g
:=e_T(p/X)~
\int_{[F^g]^{vir}}
{i_{F^g}^*\beta\over e_T(F^g/M_{g,1}(d,X))},~~~~
\beta:={\rho^*b_T(V_D)\over e_G(F_{O,D}/M_D)}.
$$

\corollary{For $g>0$,
$$
Lk_g=Lk_0~\times\int_{\bar M_{g,1}}
{b_T(\cH^*\otimes V|_q)~\Lambda_g(-\alpha)~e_T(\cH^*\otimes T_qX)\over
(-{\lambda\over\delta}+c_1(\cL))}.
$$
}
\proof
Restricting the bundle $\rho^*V_D$ on
$M_{g,1}(d,X)$ to the component $F^g$, we get
$$
i_{F^g}^*\rho^* V^g_d
=[H^1(C,f^*V)]
=[H^1(C_0,f_0^*V)]\oplus [H^1(C_1,f_1^*V)].
$$
So, we have
$$
i_{F^g}^*\rho^* b_T(V_D)
=i^*_{F^0}\rho^*b_T(V_{0,0;d})~b_T(\cH^*\otimes V|_q).
$$
Again, the first factor is constant on $F^g$.
Note that $F^g$  consists of orbifold points
of order $\delta$. Thus
the integral $\int_{[F^g]^{vir}}$
can be written as ${1\over\delta}\int_{\bar M_{g,1}}$.
Now applying the preceding theorem with
$\omega=\rho^* b_T(V_D)$
yields the desired result. $\Box$

In the special case $b_T=e_T$,
the linking values become
$$
Lk_g=Lk_0~\times\int_{\bar M_{g,1}}
{\prod_i \Lambda_g(-\xi_i)~\prod_j\Lambda_g(-\lambda_j)~\Lambda_g(-\alpha)\over
(-{\lambda\over\delta}+c_1(\cL))}
$$
where the $\xi_i$ and $\lambda_j$ are the weights on
the isotropic representations
$V|_q$ and $T_qX$ respectively. These integral
are nothing but Hodge integrals on $\bar M_{g,1}$.
Their values have been fully determined in \FP.

Fix a concave bundle $V$
and multiplicative class $b_T$. Consider the associated
Euler series $A(t)=e^{-H\cdot t/\alpha}\sum A_D~\nu^g~e^{d\cdot t}$
with coefficients
$$
A_D=e_*\left(
{\rho^* b_T(V_D)\cap LT_{g,1}(d,X)\over e_G(F_{O,D}/M_D)}
\right).
$$
By Lemma \NormalBundle, we see that
$$
deg_\alpha~A_D\leq -2+g.
$$

\theorem{Consider the gluing sequence $b_D=\Lambda_g(\alpha)$
and suppose that $c_1(X)>0$. Then
for a given $g$, the $A_D$ can
be reconstructed from the linking values $Lk_g$ and from finitely many degrees $d$.}
\proof
Recall that the homology class $LT_{g,1}(d,X)$ has dimension
$s=exp.dim~M_{g,1}(d,X)=(1-g)(dim~X-3)+\bra c_1(X),d\ket+1$.
Let $c=c_1(L_D)$. Then $c^k\cap LT_{g,1}(d,X)$ is of
dimension $s-k$, and so $e_*(c^k\cap LT_{g,1}(d,X))$ lies in the
group $A^T_{s-k}(X)$. But this group is zero unless $s-k\leq dim~X$.
The last condition means that
$$-k+2g\leq-\bra c_1(X),d\ket+g(dim~X-1)+2.$$
For given $g$, the right hand side is negative for all
but finitely many $d$. Suppose that $A_D$ are known for those
finitely many $d$.
Now
$$A_D=\sum_{k\geq0}\alpha^{-k-2}e_*\left(
\Lambda(\alpha)\Lambda(-\alpha)c^k\cap LT_{g,1}(d,X)\right)
=\sum_{k\geq0}(-1)^g\alpha^{-k-2+2g}e_*(
c^k\cap LT_{g,1}(d,X)).
$$
So each of the unknown $A_D$ has order $\alpha^{-2+p}$.
where $p<0$. By Theorem 4.3 \LLYIII, these $A_D$
can be reconstructed from the linking values.
$\Box$

The same argument shows that if $\{b_D\}$ is a given gluing sequence
with the property that for a given $g$, the number
$$
-\bra c_1(X),d\ket+g(dim~X-2)+2+deg~b_D
$$
is negative for all but finitely many $d$, then
the theorem above holds for this gluing sequence.

\subsec{The collapsing lemma}

\def\Pf{{\bf P}^n}
\def\Po{{\bf P}^1}
\def\pO{\Po}
\def\po{\Po}
\def\pf{\Pf}
\def\Pn{\Pf}
\def\fM{{\bf M}}
\def\cO{{\cal O}}
\def\dual{^{\vee}}
\def\lra{\longrightarrow}
\def\cL{{\cal L}}
\def\cF{{\cal F}}
\def\cA{{\cal A}}
\def\cE{{\cal E}}
\def\cR{{\cal R}}
\def\cS{{\cal S}}
\def\cV{{\cal V}}
\def\cK{{\cal K}}
\def\cL{{\cal L}}
\def\ex{^{\rm ex}}
\def\cM{{\cal M}}
\def\det{{\rm det}\,}
\def\mapright#1{\,\smash{\mathop{\lra}\limits^{#1}}\,}

\def\sub{\subset}
\def\cX{{\cal X}}
\def\sta{^{\ast}}
\def\mor{{\rm Mor}\,}
\def\mh{\!:\!}

\def\cM{{\cal M}}
\def\upmo{^{-1}}

\def\CC{{\bf C}}
\def\opno{\cO_{\Pn}(1)}

Let $X=\Po\times\Pn$ and let $p_1,p_2$ be the first and the second
projection of
$\pO\times\Pn$. We let $M_g(d,X)$ be the moduli space (stack) of stable
morphisms from genus $g$ curves to $X$ of bi-degree $(1,d)$.
The case $g=0$ was treated in \LLYI. Here we will prove a similar
lemma in case
$g\geq 1$. Note that there are no degree 1 maps from positive genus
smooth curves to $\Po$. Thus the domain of any stable morphism $f\mh C\to
X$ in $M_g(d,X)$
must have a distinguished irreducible component $C_0\cong\Po$ with
$$p_1\circ f|_{C_0}: C_0\mapright{\sim} \Po
$$
and all other components mapping to points via $p_1\circ f$. Let $d_0$ be
the degree of
$p_1\circ f|_{C_0}$. Use the collapsing map $M_0(d_0,X)\to N_{d_0}$, which
depend on a
choice of basis of $H^0(\opno)$, we obtain
$(n+1)$ sections
$$[\phi_0,\cdots,\phi_n]\in H^0(\cO_{\Po}(d_0))^{\oplus(n+1)}/\CC^\times.
$$
Let $C_1,\cdots,C_k$ be other
irreducible components of $C$ and let $z_i\in\Po$ be $f(C_i)$ and $d_i$ be
the degree of
$f\sta p_2\sta\cO_{\Pn}(1)$ over $C_i$.
Note $d=\sum_{i=0}^k d_i$. Then using imbedding of sheaves
$$\cO_{\Po}(d_0)\to\cO_{\Po}(\sum_{i=0}^k d_iz_i)\cong \cO_{\Po}(d)
$$
we obtain $(n+1)$-tuple of sections
$$[\tilde\phi_0,\cdots,\tilde\phi_n]\in H^0(\opno)^{\oplus(n+1)}/\CC^\times,
$$
which will be a point in $N_d$. This defines a correspondence
$$\tilde\varphi: M_g(d,X)\lra N_d.
$$

\lemma{
The correspondence $\tilde\varphi$ is induced by a
morphism $\varphi:\ M_g(d,X)\rightarrow N_d$.
Moreover $\varphi$ is equivariant with respect to the
induced action of $\CC^{\times}\times T$.
}
\proof
The following proof is given by J. Li.
Let $\cS$ be the category of all schemes of finite type (over $\bf C$)
and let $\cF: \cS\lra ({\rm Set})$
be the the contra-variant functor associating to each $S\in \cS$ the
set of families of stable morphisms
$$F: \cX\lra \po\times\pf\times S
$$
over $S$ of bi-degree $(1,d)$ of arithmetic
genus g, modulo the obvious equivalence relation. Note that
$\cF$ is represented by the moduli stack $M_g(d,X)$. Hence to define
$\varphi$
it suffices to define a transformation
$$\Psi: \cF\lra \mor(-,N_d).
$$
We now define such a transformation. Let $S\in\cS$ and let
$\xi\in \cF(S)$ be represented by
$(\cX,F)$. We let $p_i$ be the
composite of $F$ with the $i$-th projection of
$\po\times\pf\times S$ and let $p_{ij}$ be the composite of $F$ with
the projection from $\po\times\pf\times S$ to the product of its
$i$-th and $j$-th components.
We consider the locally free sheaves
$p_2\sta \cO_{\pf}(k)$, $k=0$ or $1$,
of $\cO_{\cX}$-modules and its direct image complex
$$\cL_{\xi}(k)=R^{\bullet}p_{13\ast} p_2\sta\cO_{\pf}(k).
$$
We claim that $\cL_{\xi}(k)$, which is a complex of
$\cO_{\Po\times S}$-modules, is quasi-isomorphic to a perfect complex.
Since this is a local problem, we can assume $S$ is affine.
We pick a sufficiently relative-ample line bundle $H$ on $\cX/S$ so that
$p_2\sta \cO_{\pf}(-k)$ is a quotient sheaf of
$$\cV_1=p_{3}\sta p_{3\ast}(\cO(H)\otimes p_2\sta\cO_{\Pn}(-k))\otimes
\cO(-H).
$$
Let $\cV_2$ be the kernel of $\cV_1\to p_2\sta\cO_{\Pn}(-k)$. (Here $\cV_i$
are
implicitly depending on $k$.)
Since $H$ is sufficiently ample, both $\cV_1$ and $\cV_2$ are
locally free. Hence we have a short exact sequence of locally free sheaves
of $\cO_{\cX}$-modules
$$0\lra p_2\sta\cO_{\Pn}(k) \lra \cV_1\dual\lra\cV_2\dual\lra 0.
$$
We then apply $R^{\bullet}p_{13\ast}$ to this exact sequence,
$$0\lra p_{13\ast}p_2\sta\cO_{\pf}(k)\lra
p_{13\ast}\cV_1\dual \lra p_{13\ast}\cV_2\dual\lra
R^1p_{13\ast} p_2\sta\cO_{\Pf}(k)\lra 0.
$$
Here all other terms vanish
because $H$ is sufficiently ample and fibers of $p_{13}$ have dimension at
most one.
As argued in \LLYI~ both $p_{13\ast}\cV_1\dual$ and
$p_{13\ast}\cV_2\dual$
are flat over $S$. Because $\Po\times S\to S$ is smooth,
\eqn\Two{
p_{13\ast}\cV_1\dual \lra p_{13\ast}\cV_2\dual,
}
and hence $\cL_{\xi}(k)$, is quasi-isomorphic to a
perfect complex.

The complex $\cL_{\xi}(k)$
satisfies the following
base change property: let $\rho: T\to S$ be any base change and
let $\rho\sta(\xi)\in\cF(T)$ be the pull back of $\xi$. Then there is a
canonical isomorphism of complex of sheaves of $\cO_T$-modules
$$\cL_{\rho\sta(\xi)}(k)\cong ({\bf 1}_{\po}\times \rho)\sta
\cL_{\xi}(k).
$$

Since $\cL_{\xi}(k)$ is quasi-isomorphic to a perfect complex, we can
define the determinant line bundle$^4$\footnote{}{$^4$~~All materials concerning
determinant line bundle are taken from \KM.} of $\cL_{\xi}(k)$, denoted by
$\det\cL_{\xi}(k)$.
It is an invertible sheaf of $\cO_{\po\times S}$-modules.
Using the Riemann-Roch theorem, one computes that the
degree of $\det(\cL_{\xi}(k))$ along fibers of $\Po\times S\to S$ are
$kd-g$.
Further, because $\cL_{\xi}(k)$ has rank one, there is a canonical
homomorphism
$$
\cL_{\xi}(k)\lra\det\cL_{\xi}(k)
$$
defined away from the support of the torsions subsheaves
of $p_{13\ast}p_2\sta\cO_{\Pn}(k)$
and $R^1p_{13\ast}p_2\sta\cO_{\Pn}(k)$.
Now let $w$ be any element in $H^0(\pf,\cO_{\pf}(1))$.
Its pull back provides a canonical meromorphic section
$$\sigma_{\xi,w} \in
H^0\left(\Po\times S,\fM(\det\cL_{\xi}(1))\right).
$$
For similar reason, the section $1\in H^0(\cO_{\Pn})$
provides a canonical meromorphic section $\delta$ of $\det\cL_{\xi}(0)$.
Combined, they provide a canonical meromorphic section
$$\eta_{\xi,w}=\sigma_{\xi,w}\cdot\delta\upmo \in
H^0\left(\Po\times
S,\fM(\det\cL_{\xi}(1)\otimes\det\cL_{\xi}(0)\upmo)\right).
$$
We now show that $\eta_{\xi,w}$ extends to a regular section.
Let $s\in S$ be any closed point. We first assume that there are no
irreducible components
of $\cX_s$ that are mapped entirely to $\Po\times
w\upmo(0)\sub\Po\times\Pn$ under
$F_s$. Here $F_s\mh \cX_s\to\Po\times\Pn$ is the
restriction of $F$ to the fiber over $t\in S$.
By shrinking $S$ if necessary, we can assume all
$F_t\mh\cX_t\to\Po\times\Pn$, $t\in S$,
have this property. Then the section $w$ induces a short
exact sequence
$$0\lra \cO_{\cX}\lra p_2\sta\opno\lra \cR\lra 0
$$
and a long exact sequence
$$ \lra R^{\bullet}p_{13\ast}\cO_{\cX}\lra
R^{\bullet}p_{13\ast}p_2\sta\opno\lra
R^{\bullet}p_{13\ast}\cR\lra R^{\bullet+1}p_{13\ast}\cO_{\cX}\lra .
$$
By our assumption on $w$, $R^{1}p_{13\ast}\cR=0$.
Next, we let $\cE_1\to\cE_2$ (resp. $\cF_1\to\cF_2$) be the complex
\Two~
associated to $k=1$ (resp. $k=0$). Clearly, we have canonical commutative
diagram
$$\matrix{
0&\lra& p_{13\ast}\cO_{\cX}&\lra&\cF_1&\lra&\cF_2&\lra&
R^1p_{13\ast}\cO_{\cX}&\lra&0\cr
&&\downarrow& &\downarrow&&\downarrow&&\downarrow&& \cr
0&\lra& p_{13\ast}p_2\sta\opno&\lra&\cE_1&\lra&\cE_2&\lra&
R^1p_{13\ast}p_2\sta\opno&\lra&0
}
$$
of short exact sequences.
Let
$$\cK_1:\quad \cF_1\lra\cE_1\oplus\cF_2\lra\cE_2
$$
be the induced complex. Note the last arrow is surjective. Let $\cA_1$ be
$\cF_1$ and
$\cA_2$ be the kernel of the last arrow of the above complex. Hence $\cK_1$
is quasi-isomorphic to the complex
\eqn\Five{
\cK_2:\quad\cA_1\lra \cA_2.
}
Note that both $\cA_1$ and $\cA_2$ are $\cO_S$-flat. Hence we can define
the
determinant line bundle $\det\cK_2$.
We then have canonical isomorphisms
\eqn\Four{
\det\cK_2\cong\det\cK_1\cong
\det\cL_{\xi}(1)\upmo\otimes\det\cL_{\xi}(0).
}
Now let $t\in S$ be any closed point.
$t\in S$.
It is clear that the restriction of \Five~ to
general points of $\Po_t$ is an isomorphism. Hence $\det\cK_2\upmo$
has a canonical section over $\Po\times S$ \FriedmanMorgan. It is direct to check that
under
the isomorphism \Four~ this section is the extension of $\eta_{\xi,w}$.
Since
$\po\times S\to S$ is proper, such extension is unique.

It remains to show that $\eta_{\xi,w}$ can be extended even the assumption
on $w\upmo(0)$
does not hold. Note that in this case, we can find two sections $w_1$ and
$w_2$ in
$H^0(\opno)$ so that $w=w_1+w_2$ and that both $w_1$ and $w_2$ satisfies
the condition about
$w_1\upmo(0)$ and $w_2\upmo(0)$. Here we might need to shrink $S$ if
necessary. Then
$\eta_{\xi,w_1}$ and $\eta_{\xi,w_2}$ both can be extended to regular
sections in
$$H^0\left(\Po\times S,\det\cL_{\xi}(1)\otimes\det\cL_{\xi}(0)\upmo\right).
$$
Further, over the open subset $Z\sub \Po\times S$ where all
$R^{i}p_{13\ast}p_2\sta\cO_{\Pn}(k)$, $i, k=0, 1$, are torsion free, we
obviously have
$$\eta_{\xi,w}=\eta_{\xi,w_1}+\eta_{\xi,w_2}.
$$
Since $Z\cap \Po\times\{t\}\ne \emptyset$ for all $t\in S$, the right hand
side of the above
identity provides an extension of $\eta_{\xi,w}$. This proves that for any
$w\in H^0(\opno)$
the meromorphic sections $\eta_{\xi,w}$ extends to a regular section
$$\eta_{\xi,w}\ex\in H^0\left(\Po\times
S,\det\cL_{\xi}(1)\otimes\det\cL_{\xi}(0)\upmo\right).
$$
Again since $\Po\times S\to S$ is smooth and proper, the extension is
unique.

Now we define the morphism $S\to N_d$. Let
$\{w_0,\cdots,w_n\}$ be a basis
of $H^0(\opno)$. Then we obtain $(n+1)$ canonical regular sections
$$\eta_{\xi,w_0}\ex,\cdots,\eta_{\xi,w_n}\ex
\in H^0\left(\Po\times
S,\det\cL_{\xi}(1)\otimes\det\cL_{\xi}(0)\upmo\right).
$$
Hence, after fixing an isomorphism
$$
\det\cL_{\xi}(1)\otimes\det\cL_{\xi}(0)\upmo
\cong \pi_S\sta \cM\otimes\pi_{\po}\sta\cO_{\po}(d)
$$
for some invertible sheaf $\cM$ of $\cO_S$-modules,
we obtain $(n+1)$ canonical sections of
$$
H^0(\cO_{\po}(d)) \otimes_{{\bf C}}\cM
$$
which defines a morphism
$$S\lra H^0(\opno)^{\oplus(n+1)}/\CC\sta.
$$
Since the image is always away from $0$, it defines a morphism
$$S\lra N_d.
$$
It is routine to check that this construction satisfies the base change
property,
and hence defines a morphism $M_g(d,X)\to N_d$, as desired.

To check that this morphism gives rise to the correspondence mentioned
before, it
suffices to check the case where $S$ is a closed point.
In this case one sees immediately that
the complex $\cL_{\xi}(1)-\cL_{\xi}(0)$ has locally free part isomorphic to
$\cO_{\Po}(d_0)-\cO_{\Po}$, and has torsion part supported at $z_i$ of
length $d_i$. Further, a direct check shows that the sections
$\sigma_{\xi,w_i}\cdot\delta\upmo$ are a non-zero constant multiple
of the sections $\tilde\phi_i$ mentioned in the definition of the
correspondence.
This shows that the morphism defines the correspondence constructed.
The equivariant property of this morphism again follows
from the base change property of this construction. This completes the
proof
of the Collapsing Lemma. $\Box$

\footatend\vfill\supereject\immediate\closeout\rfile\writestoppt
\baselineskip=14pt\centerline{{\bf References}}\bigskip{\frenchspacing%
\parindent=20pt\escapechar=` \input refs.tmp\vfill\eject}\nonfrenchspacing
\end